\documentclass[righttag,12pt]{amsart}
\usepackage{amssymb}
\usepackage{euscript}

\newtheorem{prop}{Proposition}[section]
\newtheorem{theorem}[prop]{Theorem}
\newtheorem{cor}[prop]{Corollary}
\newtheorem{lem}[prop]{Lemma}
 
\numberwithin{equation}{section}

\theoremstyle{definition}
\newtheorem{defn}[prop]{Definition}

\theoremstyle{remark}
\newtheorem{rem}[prop]{Remark}
\newcommand{\bbN}{{\mathbb{N}}}
 \newcommand{\bbR}{{\mathbb{R}}}
\newcommand{\bbZ}{{\mathbb{Z}}}

\newcommand{\al}{\alpha}
\newcommand{\be}{\beta}

\newcommand{\e}{\varepsilon}

\newcommand{\la}{\lambda}

\newcommand{\f}{\varphi}

\newcommand{\De}{\Delta}

\renewcommand{\span}{\operatorname{span}}
\newcommand{\supp}{\operatorname{supp}}

\newcommand{\card}{\operatorname{card}}

\newcommand{\lb}{\label}

\newcommand{\lra}{\longrightarrow}

\newcommand{\wtw}{if and only if }

\begin{document}

\title{Contractive projections in Orlicz sequence spaces}

\author{Beata Randrianantoanina$^*$ $\!^\dag$}\thanks{$^*$Participant, NSF Workshop
in Linear Analysis and Probability, Texas A\&M University}\thanks{$\dag$Partially
funded by a NSF-AWM Travel Grant}

\address{Department of Mathematics and Statistics \\ Miami University
\\Oxford, OH 45056}

\email{randrib@muohio.edu}



\begin{abstract}
We characterize  norm one complemented subspaces of    Orlicz sequence spaces $\ell_M$
equipped with either Luxemburg or Orlicz norm, provided
 that the Orlicz function $M$ is  sufficiently smooth   and
  sufficiently different from the square
function. This paper concentrates on the more difficult real case, the complex case
follows from previously known results.
\end{abstract}


\maketitle

\section{Introduction}

One of the main topics in the study of Banach spaces has been, since the inception of
the field, the study of  projections and complemented subspaces. Naturally, one of the
most important topics of the isometric Banach space theory is the study of contractive
projections (i.e. projections of norm one) and 1-complemented subspaces (i.e. ranges
of norm one projections). They were also investigated from the approximation theory
 point of view, as part of of the study of minimal projections, i.e. projections onto
 the given subspace with the smallest possible norm, for an overview of this line of
 research see \cite{ChP70,OL90}. Contractive projections are also closely related to
the metric projections or nearest point mappings, and are a natural extension of the
notion of orthogonal projections from the Hilbert space setting to general Banach
spaces. We refer the reader to the survey \cite{survey} for an outline of the
development and applications of this theory. Here we just indicate some main facts
putting the results of the present paper in context.

It is well known that in Lebesgue spaces $L_p$ and $\ell_p$, $1\le p<\infty$, a
subspace $Y$ is 1-complemented \wtw $Y$ is isometrically isomorphic to an $L_p-$space
of appropriate dimension (see \cite{Ando,D}). This is no longer the case for other
spaces. Lindberg \cite{L73} demonstrated that there exist classes of Orlicz sequence
spaces $\ell_M$  containing 1-complemented subspaces which are not even isomorphic to
$\ell_M$. In fact, he showed that for all $1<a\le b<\infty$, there exists a reflexive
Orlicz sequence space $\ell_M$ so that for all $p\in[a,b]$ there is a contractive
projection from $\ell_M$ onto a subspace isomorphic to $\ell_p$. This implies in
particular that Orlicz sequence spaces can have continuum isomorphic types of
1-complemented subspaces and thus any attempt for a geometric characterization of
1-complemented subspaces seemed hopeless.

On the other hand, 1-complemented subspacesof $\ell_p$ are also characterized as
subspaces which are spanned by a family of mutually disjoint elements of $\ell_p$ (see
\cite{LT1,BCh}). Moreover all known examples of
 1-complemented
subspaces in symmetric Banach spaces with 1-unconditional bases, and
  sufficiently different from Hilbert spaces,  are spanned by a family of
mutually disjoint vectors. (Note here that, since in Hilbert spaces every subspace is
1-complemented, it is both natural and necessary to include in this context some kind
of an assumption about the space being different from Hilbert space.) In particular,
the above described example of Lindberg of 1-complemented subspaces of Orlicz sequence
spaces which were pathological in the isomorphic sense, are not pathological in the
sense that they are spanned by   mutually disjoint vectors and the norm one projection
is the most natural averaging projection. It was shown in \cite{complex} that indeed
every 1-complemented subspace $Y$ in any complex Banach space $X$ with a
1-unconditional basis (not necessarily symmetric) which does not contain a
1-complemented isometric copy of a 2-dimensional Hilbert space $\ell_2^2$, has to be
spanned by a family of disjointly supported elements of $X$ and the norm one
projection from $X$ onto $Y$ has to be the  averaging projection. In particular, this
holds in complex Orlicz sequence spaces $\ell_M$ equipped with either the Luxemburg or
the Orlicz norm when $M$ is sufficiently different from the square function (cf.
Remark~\ref{lastr}).

In the real case this statement in its full generality is false (cf. \cite{complex}).
For real spaces  we only had the following much less satisfactory result describing
  special 1-complemented subspaces of finite codimension in Orlicz sequence spaces
$\ell_M$.

\begin{theorem}\lb{first}
\cite[Theorem~7]{real}
Let $M$ be an Orlicz function such that $M(t)>0$ for all $t>0$ and
$M$ is not similar to
$t^2$ (i.e. there do not exist constants $C, t_0 > 0$ so that
$M(t)=Ct^2$ for all $t<t_0$).
Let $\ell_M$ be the Orlicz space  equipped with either the Luxemburg or
the Orlicz
norm and $F \subset \ell_M$ be a subspace of finite codimension.
  If $F$ contains at least one basis vector and $F$ is
$1$-complemented in $\ell_M$ then  $F$ is spanned by a family of
disjointly supported vectors.
\end{theorem}

In the present paper we prove a much stronger result -- we eliminate the assumption
that the subspace should be of finite codimension. Namely we show that when $M$ is a
sufficiently smooth Orlicz function which satisfies condition $\De_2$ and is
sufficiently different from the square function, then every 1-complemented subspace of
the real Orlicz space $\ell_M$ is
 spanned by a family of mutually disjoint vectors and every
norm one projection in $\ell_M$ is an averaging projection (see Theorem~\ref{orlicz}
 and Corollary~\ref{cororlicz}). This result is valid in Orlicz spaces
   equipped with either the Luxemburg or the  Orlicz
norm.

Our method of proof is different from that of \cite{real}, it relies on new results
characterizing averaging projections through properties related to and generalizing
disjointness preserving operators \cite{aver}.

Recently, Jamison, Kami\'nska and Lewicki \cite{JKLe} obtained (using different
techniques) a generalization of Theorem~\ref{first} in another direction -- they
characterized 1-complemented subspaces of finite codimension in sufficiently smooth
Musielak-Orlicz sequence spaces, whose Orlicz function   is sufficiently different
from the square function.

We follow standard definitions and notations as may be found in \cite{Kr-Rut,LT1}

\section{Preliminary definitions}

Orlicz spaces are one of the most natural generalizations of
classical spaces $L_p$. They were first considered by Orlicz in
1930s.  Since then they were extensively studied by many authors,
see, for example the monographs \cite{Kr-Rut,RaoRen,Chen}.
 Below we recall the basic definitions and facts about Orlicz spaces that will be important
for the
present paper.

\begin{defn} We say that a function $M:  \bbR \lra [0, \infty)$ is an
{\it Orlicz
function} if $M$ is even, continuous, convex, $M(0)=0$, $M(1)=1$,
$\lim_{u\to0}M(u)/u=0$ and $\lim_{u\to\infty}M(u)/u=\infty$.
\end{defn}

Note that since the Orlicz function $M$ is convex, it has the
right derivative $M'$.  Let
 $q$ be the right inverse of $M'$.  Then we call
 $$M^*(v)=\int_0^{|v|}q(s)ds $$
 the {\it
 complementary function of } $M$.
Function $M^*$ is also an Orlicz function.
\begin{defn}
We say that the Orlicz function $M$ {\it satisfies the $\Delta_2$ condition
near zero} ($M\in \Delta_2$) if there exist constants $k>0$ and
$u_0 \geq 0$ such that for all $u$ with $|u| \leq u_0$ $$ M(2u)
\leq k M(u).  $$
\end{defn}
Note that $M\in \Delta_2$ does not imply that $M^*\in \Delta_2$.

The Orlicz function $M$ generates the {\it modular} defined for
scalar sequences $x=(x_j)_{j\in\bbN}$ by:
$$ \rho_M(x) =
\sum_{j=1}^\infty M(x_j) .  $$

The {\it Orlicz sequence space} $\ell_M$ is the space of
sequences $x$ such that there
 exists $\la>0$ with $\rho_M(\la x)<\infty$.
 If $M\in \Delta_2$ then $\ell_M= \{x\ :\
 \rho_M(\la x)<\infty \ \text{for all} \ \la\in\bbR\}.$
 The Orlicz {sequence space}
 $\ell_M$ is usually equipped with one of the two following
 equivalent norms:
 \begin{itemize}
 \item[(1)] the {\it Luxemburg norm} defined by:
 $$ \Vert x \Vert_M =
 \inf\{\lambda:  \rho_M\left(\frac{x}{\la}\right) \leq 1\}, $$
 \item[(2)] the {\it Orlicz
 norm} defined by:
 $$ \Vert x \Vert_M^O = \sup\{\sum_{j=1}^\infty x_jy_j \ :\
 \rho_{M^*}(y) \leq 1\}.  $$
 \end{itemize}

If $M\in \Delta_2$ then these norms are dual to each other in the
following sense:
\begin{eqnarray}\lb{dual1} (\ell_M,\|\cdot\|_M)^*&=&(\ell_{M^*},\|\cdot\|_{M^*}^O),\\
\lb{dual2} (\ell_M,\|\cdot\|_M^O)^*&=&(\ell_{M^*},\|\cdot\|_{M^*}).
\end{eqnarray}

We say that two Orlicz functions $M_1$ and $M_2$ are {\it
equivalent} if there exist $u_0> 0,\ k, l> 0$ such that for all
$u$ with $|u| \le u_0$ $$ M_2(ku)\le M_1(u)\le M(lu).  $$

This condition is of importance since Orlicz spaces $\ell_{M_1},
\ell_{M_2}$ are isomorphic if and only if the Orlicz functions
$M_1$, $M_2$ are equivalent.  We note that if an Orlicz function
$M$ satisfies the condition $\Delta_2$ near zero then every Orlicz
function $M_1$ equivalent to $M$ also satisfies the condition
$\Delta_2$ near zero.

  Krasnoselskii and Rutickii proved the following characterization of the
$\Delta_2$-condition in terms of the right derivative $M'$ of $M$.

\begin{prop} \lb{KRd2}\cite[Theorem 4.1]{Kr-Rut}
A necessary and sufficient condition
that the Orlicz function $M(u)$ satisfy the $\Delta_2$-condition
near zero is that there exist constants $\alpha$ and $u_0 \geq0$
such that, for $0\le u \leq u_0$
\begin{equation}\lb{delta2} \frac{u M'(u)}{M(u)} < \alpha, \end{equation}
where $M'$
denotes the right derivative of $M$.

Moreover, if \eqref{delta2} is satisfied then $M(2u)\le2^\al
M(u)$ for $0\le u \leq u_0/2$.  \end{prop}

In \cite{isoorlicz} we introduced another condition which on one
hand is very similar to \eqref{delta2}, but on the other hand is
in its nature of ``smoothness type'', as we explain below.

\begin{defn}\lb{defdelta2+} Assume that the Orlicz function $M $ is
twice differentiable
and that $M $ satisfies the $\Delta_2-$condition near zero.  We
say that $M$ {\it satisfies condition $\Delta_{2+}$ near zero} if
there exist constants $\be> 0 $ and $u_0\ge 0$ such that for all
$ u\le u_0$
\begin{equation}\lb{delta2+} \frac{u
M''(u)}{M'(u)}<\beta.
\end{equation}
\end{defn}

Condition $\Delta_{2+}$ is of ``smoothness type'' in the
following sense:
 \begin{itemize} \item[(i)] for every function $M$ which satisfies
 condition $\Delta_{2}$
 there exists an equivalent Orlicz function $M_1$ which does satisfy
 $\Delta_{2+}$;
 However, we do not know whether for every $\e>0$ it is possible to
 choose $M_1$ so that
 it is $(1+\e)-$equivalent with $M$, \item[(ii)] for every Orlicz
 function $M$ which
 satisfies $\Delta_{2+}$ there exists an equivalent (even up to an
 arbitrary $\e>0$)
 Orlicz function $M_1$ which does not satisfy $\Delta_{2+}$.
 \end{itemize}

We say that a Banach space $X$ is {\it smooth} if every element
$x\in X$ has a unique norming functional $x^*\in X^*$, i.e.  the
functional with the property that
$\|x^*\|_{X^*}^2=\|x\|_X^2=x^*(x)$ is determined uniquely for
every $x\in X$.

If $M\in\De_2$ then an Orlicz space $\ell_M$ is smooth whenever
$M$ is differentiable everywhere.

It is well known (see e.g.  \cite{Chen}) that any Orlicz function
$M$ can be ``smoothed out'', that is for any $M$ there exists an
equivalent Orlicz function $M_1 $ such that $M_1$ is twice
differentiable everywhere, $M_1''$ is continuous on $\bbR$ and
$M_1'' (u)> 0$ for all $u>0$. We recall here that a class of
functions whose second derivative exists and is continuous on
$\bbR$ is denoted by $C^2$. Thus $M_1$ above belongs to $C^2$.
Moreover, given any $\e> 0$ it is possible to choose $M_1$ so that
$\ell_M$ and $\ell_{M_1}$ are $(1+\e)-$isomorphic to each other
\cite{Chen}.

Maleev and Troyanski \cite{MT91} considered a stronger notion of
smoothness in Orlicz spaces which guarantees the
differentiability of the norm.  We recall the relevant
definitions and results.

\begin{defn} \cite{MO60} (cf.  \cite[p.  143]{LT1})
To every Orlicz function M we associate
the following {\it Matuszewska-Orlicz index}:
 $$
\alpha^{0}_{M}=\sup\{p : \sup\{\frac{M(\la t)}{t^{p}M(\la)} :
\la,t \in(0,1]\}<\infty\}.  $$
\end{defn}
\begin{defn} \cite{MT91}
We say that an Orlicz function $M$ belongs to the class $AC^{k}$
at zero if
\begin{item} $(i)$ \ $\al^{0}_{M}>k$,\\ $(ii)$ \
$M^{(k)}$ is absolutely continuous in every finite interval,\\
$(iii)$ \ $t^{k+1}|M^{(k+1)}(t)|\leq c M(ct)$ a.e. in
$[0,\infty)$ for some $c>0$.
\end{item}
\end{defn}


\begin{defn} Let $X,Y$ be Banach spaces.  The function $\f:X \lra Y$
is said to be {\it
$k-$times differentiable at} $\f \in X$ if for every $j, \ 1 \leq
j \leq k,$ there exists a continuous symmetric $j-$linear form
$T_{f}^{j}:X \times \cdots \times X = X^{(j)} \lra Y$ so that:
\begin{equation}\lb{diff} \f(f + \al g) = \f(f) + \sum^{k}_{j=1}
\al^{j} T^{j}_{f}(y,\dots,y)+\sigma_{f} (|\al|^{k})
\end{equation} uniformly on $g$ from the unit sphere $S(X)$ of
$X$.

For an open set $V\subset X,\ \f \in F^{k}(V,Y)$ means $\f$ is
$k-$times differentiable at every point of $V$.  If \eqref{diff}
is fulfilled uniformly on $f$ over a set $W\subset V$ we shall
say that $\f$ is {\it $k-$times uniformly differentiable over $W$}
and shall write $\f \in UF^{k}(W,Y)$.  We say that $X$ is
$UF^{k}$-smooth if the norm in $X$ belongs to
$UF^{k}(S(X),\bbR)$.
\end{defn}

Maleev and Troyanski proved the following results about the
uniform smoothness of Orlicz sequence spaces $\ell_{M}$:

\begin{theorem}\cite[Theorem~6]{MT91}\lb{norm}
Let M be an Orlicz function satisfying
condition $\Delta_{2}$ at zero and such that $M \in AC^{k}$ at
zero.  Then $\ell_{M}$ equipped with the Luxemburg norm is
$UF^k$-smooth.
\end{theorem}

\begin{theorem}\cite[Corollary~10]{MT91}\lb{renorm}
Let $M$ be an
Orlicz function
satisfying condition $\Delta_{2}$ at zero.  Then for every $k \in
\bbN$ such that $k<\al^{0}_{M}$ there exists an Orlicz function
$\widetilde{M}$ equivalent to $M$ at zero so that
$\ell_{\widetilde{M}}$ (with the Luxemburg norm) is
$UF^{k}$-smooth.  (In particular $\ell_{\widetilde{M}}$ is
isomorphic to $\ell_{M}$.)
\end{theorem}

We do not know whether in Theorem~\ref{renorm} it is possible for any $\e>0$ to select
$\widetilde{M}$ so that $M$ and $\widetilde{M}$ are $(1+\e)-$equivalent.

Next we recall that a Banach lattice $X$ is called {\it strictly
monotone} if $\|x+y\|>\|x\|$ for all $x, y \ge 0, \ y\ne 0,$ in
$X$.

An Orlicz space $\ell_M$ with either  the Luxemburg or the  Orlicz norm is
strictly monotone whenever $M$ is strictly increasing on
$[0,\infty)$.

\section{Tools}

In this section we gather our main tools -- facts about
contractive projections and about disjointness in Orlicz spaces.

We will say that a projection $P$ on a  purely atomic Banach lattice $X$ is an {\it
averaging projection} if there exist mutually disjoint elements $\{u_{j}\}_{j \in J}$
in $X$ and functionals $\{u^*_{j}\}_{j \in J}$ in $X^{\ast}$ so that
$u_{j}^{\ast}(u_{k})=0$ if $j \ne k$, $u_j^{\ast}(u_{j}) = 1$ for all $j \in J$ and
for each $f \in X$
$$ Pf = \sum_{j \in J} \
u_{j}^{\ast}(f)u_{j}.  $$

First we  recall two abstract conditions that we introduced in \cite{aver} in our
study of averaging projections in purely atomic Banach lattices.

\begin{defn} \cite{aver}
Let $X$ be a Banach lattice and $P :  X \to X$ be a
linear operator on $X$.  We say that the operator $P$ is
\begin{itemize}
\item[(1)]   {\it semi band preserving} if and only if for all $f, g \in X$,
\begin{equation} \lb{D}
\supp(Pf) \cap\supp(g)=\emptyset {\text{ implies that }}
\supp(Pf) \cap\supp(Pg)=\emptyset.
 \end{equation}
 \item[(2)]   {\it semi containment preserving} if and only if for all $f, g
\in X$,
\begin{equation} \lb{C}
\supp g\subset\supp Pf {\text{
implies that }} \supp Pg\subset\supp Pf.
\end{equation}
\end{itemize}
\end{defn}

In the above statement all set relations are considered modulo
sets of measure zero.

It is clear that all averaging projections are both semi  band
preserving and
semi  containment   preserving.  In \cite{aver} we proved that in fact in
``nice'' purely atomic Banach spaces either of
semi  band or semi containment   preservation
 characterizes averaging projections among
contractive projections.  More precisely, we have:

\begin{theorem}\lb{smooth} \cite{aver}
Let $X$ be a purely atomic  strictly monotone Banach lattice and let $P:X \to X$ be a
norm one projection which is semi  band   preserving or semi  containment preserving.
Then $P$ is an averaging projection.
\end{theorem}

This theorem will be very useful for our considerations since in
\cite{isoorlicz} we obtained conditions which partially describe
disjointness and containment of supports of elements in Orlicz
spaces.  These conditions will enable us to verify that
contractive projections in Orlicz sequence spaces are
semi  band   preserving or
semi  containment   preserving.

  We note here that all theorems in \cite{isoorlicz} were formulated
  and proved for
Orlicz function spaces $L_M$, where $M$ is an Orlicz function
satisfying conditions $\De_2$ and $\De_{2+}$ near infinity.
However to adapt to the case of Orlicz sequence spaces $\ell_M$,
where $M$ is an Orlicz function satisfying conditions $\De_2$ and
$\De_{2+}$ near zero, the proofs require only very minor changes,
if any.  Thus in the following when we refer to the statements
from \cite{isoorlicz} we will formulate them using $\ell_M$
instead of $L_M$, which is more appropriate for the present paper.

  We stress that theorems in \cite{isoorlicz} are proven for Orlicz
  spaces equipped with
the Luxemburg norm, and the analogs of most of the results from
\cite{isoorlicz} are false in Orlicz spaces equipped with the
Orlicz norm.

\begin{prop} \cite[Proposition~3.1]{isoorlicz} \lb{disjoint}
Assume that $M$ is an Orlicz
function which satisfies condition $\Delta_{2+}$ and such that
$M''$ is a continuous function with $M''(0)=0$ and $M''(t) > 0$
for all $t > 0$.  Let $f, g\in \ell_M$ and $N(\alpha) =
\|f+\alpha g\|_M$.  Then
\begin{enumerate}
\item[(a)] If $f, g $
have disjoint supports, $\mu(\supp g)<\infty$ and $g$ is bounded
then $N'(0) = 0$ and $N''(\alpha) \lra 0$ as $\al\lra 0$ along a
subset of $[0,1]$ of full measure.
\item[(b)] If $N'(0) = 0$ and $N''(\alpha) \lra 0$ as $\al\lra 0$ along a subset of
$[0,1]$ of full measure then $f, g $ have disjoint supports.
\end{enumerate}
\end{prop}

\begin{prop} \cite[Proposition~4.1]{isoorlicz} \lb{Minfinity}
Assume that $M$ is an
Orlicz function which satisfies condition $\Delta_{2+}$ near zero
and such that $M''$ is a continuous function on $(0,\infty)$ with
$\lim_{t\to 0}M''(t) = \infty$.  Let $f, g \in \ell_M$ with $ f ,
g \neq 0$ and $N(\alpha) = \|f+\alpha g\|_M$.  Then
\begin{itemize}
\item[$(a)$] If $\mu(\supp g\setminus\supp f) > 0$ then $N''(\alpha) \lra \infty$, as
${\alpha\lra 0}$ along a subset of $[0,1]$ of full measure.
\item[$(b)$] If $g$ is simple and $\mu(\supp g \setminus \supp f)
= 0$ then there exists a subset $E$ of $[0,1]$ of full measure
and $C>0$ such that for all ${\alpha\in E}$
\begin{equation*}
N''(\alpha) \le C.
\end{equation*}
\end{itemize}
\end{prop}

\begin{rem} A careful reader may have noticed that
Proposition~\ref{Minfinity} above
appears slightly stronger than
\cite[Proposition~4.1]{isoorlicz}.  However the differences
between these two statements are minimal and result from a slight
simplification of the proof of \cite[Proposition~4.1]{isoorlicz}
in the case of sequence Orlicz spaces.  Also the formulation of
Proposition~\ref{Minfinity} clarifies a slight ambiguity of the
statement of \cite[Proposition~4.1]{isoorlicz}.  We leave the
details, which are easy but require cumbersome notation, to the
interested reader.
\end{rem}

Finally we recall a result from \cite{complex} which describes
the form of two dimensional 1-com\-ple\-ment\-ed subspaces of Orlicz
sequence spaces, when the two spanning elements have disjoint
supports.  (We say that a subspace is 1-complemented if it is the
range of a projection $P$ with $\|P\|=1$.)  This result will
allow us to give a very detailed description of 1-complemented
subspaces of any dimension of Orlicz sequence spaces.

\begin{theorem} \cite[Theorem~6.1]{complex} \label{block}
Let $M$ be an Orlicz function
satisfying condition $\Delta_2$ and $\ell_M$ be a (real or
complex) Orlicz sequence space equipped with either  the Luxemburg or
 the Orlicz norm and let $x,y \in \ell_M$, be disjoint norm one
elements such that $\span \{x,y \} $ is $1$-complemented in $\ell_M$. Then one of
three possibilities holds:
\begin{itemize}
\item[(1)] $\card (\supp x) < \infty$ and $|x_i| = |x_j|$ for all
$i,j \in \supp x$; or
\item[(2)] there exists $p$, $1 \leq p \leq
\infty$, such that $M(t) = Ct^p$ for all $t \leq \|x\|_\infty$;
or
\item[(3)] there exists $p$, $1 \leq p \leq \infty$, and
constants $C_1, C_2, \gamma \geq 0$ such that $C_2 t^p
\leq\nobreak M(t) \leq C_1 t^p$ for all $t \leq \|x\|_\infty$ and
such that, for all $j \in \supp x$, $$|x_j| = \gamma^{k(j)} \cdot
\|x\|_\infty$$ for some $k(j) \in \bbZ$.
\end{itemize}
\end{theorem}

In particular, it follows from Theorems \ref{block} that in
``most'' Orlicz spaces the
 only $1$-complemented disjointly supported subspaces of any dimension
 are those spanned
 by a block basis with constant coefficients of some permutation
 of the original basis.

\section{Main results}

We start from a lemma which will allow us to apply
Propositions~\ref{disjoint} and
 \ref{Minfinity} to study whether contractive
 projections in Orlicz
 sequence spaces are semi  band   preserving or
semi  containment   preserving.

\begin{lem}\lb{convex} Suppose that $\f, \psi :  \bbR \to [0,\infty)$ are convex
functions, differentiable everywhere and such that $\f(0) =
\psi(0), \f(\al) \le \psi(\al)$ for all $\al \in \bbR$.\newline
$(i)$ \ Then $\psi{'}(0) = \f{'}(0) $.\newline $(ii)$ \ If
$\f{''}(0),\psi{''}(0)$ exist and $\psi{''}(0) = 0$, then
$\f{''}(0) = 0$.\newline $(iii)$ Suppose that $\f'$ and $\psi'$
are absolutely continuous on $[0,1]$. Then, if $\f''(\al) \lra
\infty$ as $\al \lra 0$ along a subset of $[0,1]$ of full
measure, then for every $C>0$
$$ \mu (\{\al \in [0,1] :
\psi''(\al) \ exists \ and \ \psi''(\al)\leq C\})<1.  $$
\end{lem}

\begin{proof} To prove $(i)$ observe that, since $\f(0)=\psi(0)$, we
have for all $\al
\in \bbR$
$$ \f(\al)-\f(0) \le \psi(\al)-\psi(0).  $$
Thus for
$\al > 0$
\begin{equation}\lb{L1}
\frac{\f(\al)-\f(0)}{\al} \ \le \ \frac{\psi(\al)-\psi(0)}{\al} ,
\end{equation} and for $\al < 0$ \begin{equation}\lb{L2} \frac{\f(\al)-\f(0)}{\al} \ \ge
\ \frac{\psi(\al)-\psi(0)}{\al} .
\end{equation}

Since $\f{'}(0)$ and $\psi{'}(0)$ exist we have, by \eqref{L1},
$$ \f{'}(0)= \lim_{\al \to 0^{+}} \ \frac{\f(\al)-\f(0)}{\al} \
\le \ \lim_{\al \to 0^{+}} \ \frac{\psi(\al)-\psi(0)}{\al} =
\psi{'}(0), $$
and, by \eqref{L2},
$$ \f{'}(0)= \lim_{\al \to
0^{-}} \ \frac{\f(\al)-\f(0)}{\al} \ \ge \ \lim_{\al \to 0^{-}} \
\frac{\psi(\al)-\psi(0)}{\al} = \psi{'}(0), $$
Thus $\f{'}(0)=
\psi{'}(0)$ and $(i)$ is proved.

To prove $(ii)$, we consider the set $A = \{\al > 0 :
\f{'}(\al)=\psi{'}(\al)\}$.

If $\inf \{\al \in A\}=0$, then there exists a sequence
$\{\al_{n}\}_{n=1}^{\infty} \subset A$ so that $\lim_{n \to
\infty} \al_{n} = 0$.  Since $\f{''}(0)$ and $\psi{''}(0)$ exist,
and by $(i)$, we obtain:
$$ \f{''}(0)=\lim_{n \to \infty} \
\frac{\f{'}(\al_{n})-\f{'}(0)}{\al_{n}} \ = \lim_{n \to \infty} \
\frac{\psi{'}(\al_{n})-\psi{'}(0)}{\al_{n}} \ = \psi{''}(0)=0, $$
So $(ii)$ is proved.

If $\inf\{\al \in A\} > 0$ (this includes the case that $A =
\emptyset$ and then we say $\inf \{\al \in A\} = \infty > 0$),
then there exists $\e, \ 0 < \e<\inf\{\al \in A\}$ so that
$\f{'}(\al) \ne \psi{'}(\al)$ for all $\al \in (0,\e)$.

Let $h = \psi - \f$.  Then $h(\al) \ge 0$ for all $\al \in \bbR,
\ h(0)=0$ and $h{'}(\al) \ne 0$ for all $\al \in (0,\e)$.  Since
$h{'}$ satisfies the Darboux property, we get either:
\begin{equation}\lb{L3} h{'}(\al)>0 \ \ {\text{for all}} \ \al
\in (0,\e),
\end{equation} or
\begin{equation}\lb{L4} h{'}(\al)<0 \ \ {\text{for all}} \ \al \in
(0,\e).  \end{equation}
But $h(0)=0$ and $h(\e) \ge 0$, so by the
Mean Value Theorem there exists $\al_{0} \in (0,\e)$ so that
$$
h{'}(\al_{0}) = \frac{h(\e)}{\e} \ge 0.
$$
Thus \eqref{L3} has
to hold.  This implies that, since $h{''}(0)$ exists, $h{''}(0)
\ge 0$.  This means:
$$
0=\psi{''}(0) \ge \f{''}(0).
$$
Since
$\f$ is convex, we also get
$$ \f{''}(0) \ge 0.  $$
Thus $\f{''}(0)=0$ and $(ii)$ is proved.

To prove $(iii)$ we denote by $E_{1} = \{\al \in [0,1] :
\f''(\al) \ \text{exists}\}$.

Since $\f$ and $\psi$ are convex, $\mu(E_{1})=1$.  Without loss
of generality we can also assume that
\begin{equation}\lb{lim}
\f''(\al) \lra \infty \ \text{as} \ \al \lra 0 \ \text{and} \ \al
\in E_{1}.
\end{equation}

Suppose, for contradiction, that there exists $C>0$ so that the
set
\begin{equation*} E_{2} = \{\al \in [0,1] :  \psi''(\al)\
\text{exists and}\ \psi''(\al)\leq C\}
\end{equation*}
has full measure.  Let $E = E_{1}\cap E_{2}$.  By \eqref{lim} there
exists $\e>0$ so that:
\begin{equation}\lb{epsilon}
\f''(\al)>C
\ \ \text{for every $\al \in E \cap (0,\e)$}.  \end{equation}

Now consider the set $A = \{\al>0 :  \f'(\al) = \psi'(\al)\}$
similarly as we did in the proof of $(ii)$.  If $\inf\{\al \in
A\} = 0$, then there exist $\al_{1},\al_{2} \in (0,\e)$ so that
$\al_{1} \ne \al_{2}$ and
\begin{equation}\lb{equal} \begin{split}
\f'(\al_{1}) &= \psi'(\al_{1}),\\ \f'(\al_{2}) &=
\psi'(\al_{2}).
\end{split}
\end{equation}

But, since $\f'$ is absolutely continuous on $[0,1]$, and by
$\eqref{epsilon}$, we also have:
\begin{equation*} \f'(\al_{1})
- \f'(\al_{2}) = \int^{\al_{2}}_{\al_{1}} \f''(\al) d \al =
\int_{[\al_{1},\al_{2}] \cap E} \f''(\al) d \al> C
(\al_{1}-\al_{2}).
\end{equation*}

On the other hand, by the absolute continuity of $\psi'$ on
$[0,1]$ and the definition of $E_{2}$ we have:
\begin{equation*}
\psi'(\al_{1}) - \psi'(\al_{2}) = \int^{\al_{2}}_{\al_{1}}
\psi''(\al) d \al = \int_{[\al_{1},\al_{2}] \cap E} \psi''(\al) d
\al \leq C (\al_{1}-\al_{2}).
 \end{equation*}

This is a contradiction since $\eqref{equal}$ implies that $
\f'(\al_{1})-\f'(\al_{2}) = \psi'(\al_{1})-\psi'(\al_{2}) .  $

Now let us consider the case that $\inf\{\al \in A\} \ne 0$,
i.e.  $\inf\{\al \in A\} >
 0$ (this, as in $(ii)$, includes the possibility that $A = \emptyset$
 in which case we
 say that $\inf\{\al \in A\} = \infty)$.  We showed in the proof of
  $(ii)$ (cf.
 \eqref{L3}) that in this case there exists
 $\e_{1},\ 0<\e_{1}<\inf\{\al \in A\}$, so
 that
 \begin{equation}\lb{L5}
 \psi'(\al)>\f'(\al)\ \text{for all}\ \al \in (0,\e_{1}).
 \end{equation}

By $(i)$, $\f'(0) = \psi'(0)$.  Let $\al_{0} \in (0,\e) \cap
(0,\e_{1})$.  Then, similarly as in the previous case, since
$\f'$ is absolutely continuous on $[0,1]$, by $\eqref{epsilon}$,
we obtain:
 \begin{equation*}
 \f'(\al_{0})-\f'(0) =
\int^{\al_{0}}_{0} \f''(\al) d \al = \int_{[0,\al_{0}] \cap E}
\f''(\al) d \al > C \al_{0}.
 \end{equation*}

On the other hand, again by the absolute continuity of $\psi'$
and the definition of $E_{2}$:
\begin{equation*}
\psi'(\al_{0})-\psi'(0) = \int^{\al_{0}}_{0} \psi''(\al) d \al =
\int_{[0,\al_{0}] \cap E_{2}} \psi''(\al) d \al \leq C \al_{0}.
\end{equation*}

Thus
\begin{equation*}
 \psi'(\al_{0}) < \f'(\al_{0}),
\end{equation*}
 which contradicts \eqref{L5} and ends the proof
of $(iii)$.
\end{proof}

We are now ready for our main results.

\begin{theorem}\lb{DC}
Let $M$ be an Orlicz function which satisfies condition $\De_{2+}$
near zero  and let $\ell_M$ be the real Orlicz
sequence space equipped with the Luxemburg norm.  Suppose that
$P :  \ell_{M} \to \ell_{M}$ is a contractive projection.  Then
the following hold:
\begin{itemize}
\item[$(a)$] If $M\in AC^2$,
$M{''}(0)=0$ and $M{''}(t)>0$ for all $t>0$, then $P$ is
semi  band   preserving;
\item[$(b)$] If $M\in AC^1$
near zero, $M''$ is continuous on $(0,\infty)$ and $\lim_{t\to
0}M''(t) = \infty$, then $P$ is
semi  containment   preserving.
\end{itemize}
\end{theorem}

\begin{proof}
Since bounded functions with finite supports are linearly dense in
$\ell_{M}$, to show that $P$ is semi  band   preserving or
semi  containment   preserving, respectively, it is enough to verify that
\eqref{D} or \eqref{C}, resp., are satisfied with the additional
assumption that $g$ is a bounded function and $\mu(\supp g) <
\infty$.

For any functions $f, g \in\ell_M$ we define $$\psi(\al) = \|Pf +
\al g \|_{M}$$
$$\f(\al) = \|Pf + \al Pg \|_{M}$$ for all $\al \in \bbR$.  Then $\f$ and $\psi$ are
convex functions and $\psi(0) = \|Pf\| = \f(0)$.  Moreover, by
Theorem~\ref{norm}, in both cases $(a)$ and $(b)$, $\f$ and
$\psi$ are differentiable everywhere.  Since $P$ is a contractive
projection, we also get $\f(\al) \le \psi(\al)$ for all $\al \in
\bbR$.

Now to prove $(a)$ assume that $\mu(\supp g) < \infty$ and
$\supp(g)\cap\supp(Pf)=\emptyset$.  Since $M\in AC^2$, by
Theorem~\ref{norm}, $\f{''}(0)$ and $ \psi{''}(0)$ exist.  By
Proposition~\ref{disjoint}(a) we get $\psi{'}(0)=0$ and
$\psi{''}(0)=0$.  Hence by Lemma~\ref{convex}$(i)$ and $(ii)$,
$\f{'}(0)=0$ and $\f{''}(0)=0$.  Thus, by
Proposition~\ref{disjoint}(b), we get that $Pf$ and $Pg$ have
disjoint supports, which proves that $P$ is semi  band   preserving.

To prove $(b)$ assume, for contradiction, that there exist $f, g
\in\ell_M$ so that $\mu(\supp g) < \infty$, $\supp(g)\subseteq
\supp(Pf)$ and $\supp(Pg)\not\subseteq \supp(Pf)$.

Note that since $M\in AC^1$, by Theorem~\ref{norm}, functions
$\f$ and $\psi$ are differentiable everywhere, $\f', \psi'$ are
absolutely continuous and $\f'', \psi''$ exist almost
everywhere.  Further, by Proposition~\ref{Minfinity}$(b)$, there
exists a subset $E$ of $[0,1]$ of full measure and $C_0>0$ such
that for all ${\alpha\in E}$
\begin{equation}\lb{estimate}
 \psi''(\alpha) \le C_0.
 \end{equation}

On the other hand, by Proposition~\ref{Minfinity}$(b)$,
$\psi''(\alpha) \lra \infty$, as ${\alpha\lra 0}$ along a subset
of $[0,1]$ of full measure.  Hence, by Lemma~\ref{convex}$(iii)$
for every $C>0$ $$ \mu (\{\al \in [0,1] :  \psi''(\al) \
\text{exists and } \ \psi''(\al)\leq C\})<1.  $$ This contradicts
\eqref{estimate} and ends the proof of part $(b)$.  \end{proof}

As a consequence we obtain the characterization of contractive
projections in Orlicz sequence spaces.

\begin{theorem}\lb{orlicz}
Suppose that $M$ is an Orlicz function such that $M$ satisfies
condition $\De_{2+}$ near zero and one of the following two
conditions:
\begin{itemize}
\item[$(i)$] $M\in AC^2$, $M{''}(0)=0$ and $M{''}(t)>0$ for all $t>0$.
\item[$(ii)$] $M\in AC^1$ near zero, $M''$ is continuous on $(0,\infty)$
and $\lim_{t\to
0}M''(t) = \infty$.
\end{itemize}

Let $\ell_{M}$ be the real Orlicz sequence space equipped with
the Luxemburg norm and let $P $ be a contractive projection on
$\ell_{M}$.  Then $P$ is an averaging projection, i.e.  there
exist mutually disjoint elements $\{u_{j}\}_{j \in J}$ in
$\ell_{M}$ and functionals $\{u_{j}^{*}\}_{j \in J}$ in
$(\ell_{M})^{*}$ so that $u^{*}_{j}(u_{k})=0$ if $j \ne k$,
$u^{*}_{j}(u_{j})=1$ for all $j \in J$ and for each $f \in
\ell_{M}$.
\begin{equation}\lb{form P}
Pf = \ \sum_{j \in J} \ u_{j}^{*} (f) u_{j}.
\end{equation}

Moreover, one of the three possibilities holds:
 \begin{itemize}
\item[$(1)$] $\card(\supp u_{j})<\infty$ for each $j \in J$, and
$|(u_{j})_{k}|=|(u_{j})_{l}|$ for each $k,l \in \supp(u_{j}),j
\in J$.  (Here $u_{j}= \ \sum\limits_{k \in \supp u_{j}} \
(u_{j})_{k} e_{k})$; or
\item[$(2)$] there exist $p,\ 1 <
p<\infty$, and $C \in \bbR$, so that $M(t)=Ct^{p}$ for all $t \le
\sup\limits_{j \in J}\|u_{j}\|_{\infty} (\le \infty)$; or
\item[$(3)$] there exist $p,\ 1 <p< \infty$, and constants
$C_{1},C_{2},\gamma >0$, so that $C_{2}t^{p} \le M(t) \le C_{1}
t^{p}$ for all $t \le \sup\limits_{j \in J}\|u_{j}\|_{\infty}(\le
\infty),\ \|u_{j}\|_{\infty}<\infty$ for all $j \in J$, and $$
|(u_{j})_{k}| \in \{ \gamma^{m} \cdot \|u_{j}\|_{\infty} :  m \in
\bbZ \} $$ for all $j \in J$ and $k \in \supp(u_{j})$.
\end{itemize}
\end{theorem}

\begin{proof} Note first that either condition $(i)$ or $(ii)$ implies that
$\ell_M$ is smooth and that $M'$ is a strictly increasing
function on $(0,\infty)$. Thus $M$ is also strictly increasing
on $(0,\infty)$ and $\ell_M$ is strictly monotone.
Hence the fact that $P$ is an averaging projection follows
immediately from Corollary~\ref{smooth} and Proposition~\ref{DC}.

The moreover part follows directly from
\cite[Theorem~6.1]{complex} (see Theorem~\ref{block}).  Indeed,
since the elements $\{u_{j}\}_{j \in J}$ are mutually disjoint,
for any $j_1, j_2 \in J$ and any $f \in \ell_{M}$ we have
\begin{equation*}
\|u_{j_1}^{*} (f) u_{j_1}+u_{j_2}^{*} (f)
u_{j_2}\| \le \|\sum_{j \in J} \ u_{j}^{*} (f) u_{j}\|=\|Pf\|\le
\|f\|.
\end{equation*}
Thus the projection
$Q:\ell_M\lra\span\{u_{j_1}, u_{j_2}\}$ defined by $Qf =
u_{j_1}^{*} (f) u_{j_1}+u_{j_2}^{*} (f) u_{j_2}$, has $\|Q\|=1$.
Thus, by Theorem~\ref{block}, conditions $(1)$-$(3)$ in the
statement of Theorem~\ref{orlicz} are satisfied.
\end{proof}

By duality we also obtain the description of contractive
projections in real Orlicz sequence spaces equipped with the
Orlicz norm.

\begin{cor}\lb{cororlicz}
Suppose that $M$ is an Orlicz function such that $M$ satisfies
condition $\De_2$ near zero and $M^*$ satisfies condition
$\De_{2+}$ near zero and one of the following two conditions:
\begin{itemize}
\item[$(i^*)$] $M^*\in AC^2$ near zero, $M''$ is continuous on
$(0,\infty)$, $M{''}(t)>0$ for all $t>0$ and $\lim_{t\to 0}M''(t)
= \infty$.
\item[$(ii^*)$] $M^*\in AC^1$ near zero, $M\in C^2$, $M{''}(t)>0$ for all $t>0$ and
$M{''}(0)=0$
\end{itemize}

Let $\ell_{M}$ be the real Orlicz sequence space equipped with
the Orlicz norm and let $P $ be a contractive projection on
$\ell_{M}$.  Then $P$ has the form described in
Theorem~\ref{orlicz}.
\end{cor}

\begin{proof}
This follows from Theorem~\ref{orlicz} by duality. Indeed, since
$M\in \De_2$, by \eqref{dual2} we have
$(\ell_M,\|\cdot\|_M^O)^*=(\ell_{M^*},\|\cdot\|_{M^*})$ and the
dual projection $P^*$ is contractive in $\ell_{M^*}$ equipped
with the Luxemburg norm.
Further, either of the conditions $(i^*)$ or $(ii^*)$ implies
that $M^*$ is smooth, so the only thing that needs to be verified
is that condition $(i^*)$  implies that $M^*$ satisfies condition
$(i)$ and condition $(ii^*)$  implies that $M^*$ satisfies
condition $(ii)$ from Theorem~\ref{orlicz}.

For that, note that by the definition of the complementary
function $M^*$ and since in either case $(i^*)$ or $(ii^*)$,
$M''(t)> 0$ for $t> 0$, we have for all $t> 0$
$$(M^*)''(t)=\frac1{M''((M^*)'(t))}.$$

Since $M''$ and $(M^*)'$ are both continuous on $(0,\infty)$ in
either case $(i^*)$ or $(ii^*)$, we conclude that also $(M^*)''$
is continuous on $(0,\infty)$ and $(M^*)''(t)>0$ for all $t> 0$.

Moreover, since $\lim_{t\to 0} (M^*)'(t) = (M^*)'(0) = 0$, we have
in case $(i^*)$:
$$\lim_{t\to 0} (M^*)''(t) =\lim_{t\to 0} \frac1{M''((M^*)'(t))} = \lim_{s\to 0} \frac1{M''(s)} =0.$$

It is not difficult to check that this implies that
$(M^*)''(0)=0$. Therefore condition $(i)$ is implied by $(i^*)$.

Similarly, in  case $(ii^*)$ we have:
$$\lim_{t\to 0} (M^*)''(t) =\lim_{t\to 0} \frac1{M''((M^*)'(t))} = \lim_{s\to 0} \frac1{M''(s)} =\infty.$$

So condition $(ii)$ is implied by $(ii^*)$.

Hence, by Theorem~\ref{orlicz}, in  either case $(i^*)$ or $(ii^*)$,
$P^*$, and thus also $P$, have form \eqref{form P} and the conditions
$(1)-(3)$ from Theorem~\ref{orlicz} hold.
\end{proof}

\begin{rem} \lb{lastr}
We do not know whether the assumption about smoothness of $M$ is
necessary for Theorem~\ref{orlicz} and Corollary~\ref{cororlicz} to hold.
We suspect that, similarly as in the complex case, smoothness of $M$
should not be necessary.

However it is clear that some assumption about a behavior of $M''$ near zero
is necessary. Indeed in \cite[Example~3]{complex} we showed that if
$a\in(\sqrt{2/3},1)$ and
$$M_a (t) = \left\{
\begin{array}{ll}
t^2&\mbox{if }\ 0\le t\le a,\\
(1+a)t-a&\mbox{if }\ a\le t\le1,\end{array}
\right.$$
then the real or complex 4-dimensional Orlicz space $\ell_{M_a}^4$
equipped
with either the Luxemburg or the Orlicz norm contains a 2-dimensional
1-complemented isometric copy of $\ell_2^2$ which cannot be spanned by
a family of disjoint vectors from $\ell_{M_a}^4$.
It is not difficult to adjust this example so that if $a$ is any positive
number then the real or complex  Orlicz space $\ell_{M_a}$
(of infinite dimension)
 contains a 2-dimensional
1-complemented isometric copy of $\ell_2^2$ which cannot be spanned by
a family of disjoint vectors from $\ell_{M_a}$.

It would be interesting to characterize what condition on $M$ is equivalent to the
fact that $\ell_{M}$ (complex or real) does not
 contain a 2-dimensional
1-complemented isometric copy of $\ell_2^2$ (which cannot be spanned by a family of
disjoint vectors from $\ell_{M}$). Either of the conditions $(i)$, $(ii)$, $(i^*)$ or
$(ii^*)$ is clearly sufficient, but they all involve smoothness. We suspect that the
right condition is that for all $a>0$ the function $M(t)/t^2$ is not constant on the
interval $(0,a)$.
\end{rem}


\begin{thebibliography}{10}

\bibitem{Ando}
{\sc T.~Ando}, {\em Contractive projections in {$L_p-$}spaces}, Pacific J.
  Math., 17 (1966), pp.~391--405.

\bibitem{BCh}
{\sc J.~Blatter and E.~W. Cheney}, {\em Minimal projections on hyperplanes in
  sequence spaces}, Annali di Mat. Pura ed Appl., 101 (1974), pp.~215--227.

\bibitem{Chen}
{\sc S.~Chen}, {\em {Geometry of Orlicz spaces}}, Dissertationes Math., 356
  (1996), pp.~1--204.

\bibitem{ChP70}
{\sc E.~W. Cheney and K.~H. Price}, {\em Minimal projections}, in Approximation
  Theory (Proc. Sympos., Lancaster, 1969), Academic Press, London, 1970,
  pp.~261--289.

\bibitem{D}
{\sc R.~G. Douglas}, {\em Contractive projections on an ${L}_1$ space}, Pacific
  J. Math., 15 (1965), pp.~443--462.

\bibitem{JKLe}
{\sc J.~E. Jamison, A.~{Kami\'nska}, and G.~Lewicki}, {\em One complemented
  subspaces of {M}usielak-{O}rlicz sequence spaces}.
\newblock preprint.

\bibitem{Kr-Rut}
{\sc M.~A. Krasnosel'skii and Y.~B. Rutickii}, {\em Convex functions and
  {O}rlicz spaces}, P. Noordhoff LTD., Groningen, The Netherlands, 1961.

\bibitem{L73}
{\sc K.~Lindberg}, {\em On subspaces of {O}rlicz sequence spaces}, Studia
  Math., 45 (1973), pp.~119--146.

\bibitem{LT1}
{\sc J.~Lindenstrauss and L.~Tzafriri}, {\em { Classical {B}anach spaces, Vol.
  1, {S}equence spaces}}, Springer--Verlag, Berlin--Heidelberg--New York, 1978.

\bibitem{MT91}
{\sc R.~P. Maleev and S.~L. Troyanski}, {\em Smooth norms in {O}rlicz spaces},
  Canad. Math. Bull., 34 (1991), pp.~74--82.

\bibitem{MO60}
{\sc W.~Matuszewska and W.~Orlicz}, {\em On certain properties of $\varphi
  $-functions}, Bull. Acad. Polon. Sci. S\'er. Sci. Math. Astronom. Phys., 8
  (1960), pp.~439--443.

\bibitem{OL90}
{\sc W.~Odyniec~[V. P.~Odinec] and G.~Lewicki}, {\em Minimal projections in
  {B}anach spaces, {P}roblems of existence and uniqueness and their
  application}, Lecture Notes in Mathematics, 1449. Springer-Verlag, Berlin,
  1990.


\bibitem{complex}
{\sc B.~Randrianantoanina}, {\em $1$-complemented subspaces of spaces with
  $1$-unconditional bases}, Canad. J. Math., 49 (1997), pp.~1242--1264.

\bibitem{real}
{\sc B.~Randrianantoanina}, {\em One-complemented
  subspaces of real sequence spaces}, Results Math., 33 (1998), pp.~139--154.

\bibitem{isoorlicz}
{\sc B.~Randrianantoanina}, {\em Injective isometries
  in {O}rlicz spaces}, in Function spaces (Edwardsville, IL, 1998), Amer. Math.
  Soc., Providence, RI, 1999, pp.~269--287.

\bibitem{survey}
{\sc B.~Randrianantoanina}, {\em Norm one projections in {B}anach spaces},
  Taiwaneese J. Math., 5 (2001), pp.~35--95.

\bibitem{aver}
{\sc B.~Randrianantoanina}, {\em A disjointness type property of conditional
  expectation operators}.
\newblock preprint, available on the Mathematics ArXiv at
  http://front.math.ucdavis.edu/math.FA/0112181.

\bibitem{RaoRen}
{\sc M.~M. Rao and Z.~D. Ren}, {\em Theory of {O}rlicz spaces}, Marcel Dekker
  Inc., New York, 1991.

\end{thebibliography}

\def\polhk#1{\setbox0=\hbox{#1}{\ooalign{\hidewidth
  \lower1.5ex\hbox{`}\hidewidth\crcr\unhbox0}}} \def\cprime{$'$}

\end{document}